# COMBINING INFORMATION FROM INDEPENDENT SOURCES THROUGH CONFIDENCE DISTRIBUTIONS[1]

By Kesar Singh, Minge Xie and William E. Strawderman

*Rutgers University*

This paper develops new methodology, together with related theories, for combining information from independent studies through confidence distributions. A formal definition of a confidence distribution and its asymptotic counterpart (i.e., asymptotic confidence distribution) are given and illustrated in the context of combining information. Two general combination methods are developed: the first along the lines of combining $p$-values, with some notable differences in regard to optimality of Bahadur type efficiency; the second by multiplying and normalizing confidence densities. The latter approach is inspired by the common approach of multiplying likelihood functions for combining parametric information. The paper also develops adaptive combining methods, with supporting asymptotic theory which should be of practical interest. The key point of the adaptive development is that the methods attempt to combine only the correct information, downweighting or excluding studies containing little or wrong information about the true parameter of interest. The combination methodologies are illustrated in simulated and real data examples with a variety of applications.

**1. Introduction and motivations.** Point estimators, confidence intervals and $p$-values have long been fundamental tools for frequentist statisticians. Confidence distributions (CDs), which can be viewed as "distribution estimators," are often convenient devices for constructing all the above statistical procedures plus more. The basic notion of CDs traces back to the fiducial distribution of Fisher (1930); however, it can be viewed as a pure frequentist concept. Indeed, as pointed out in Schweder and Hjort (2002),

Received May 2003; revised April 2004.

[1]Supported in part by NSF Grant SES-02-41859.

the CD concept is "Neymannian interpretation of Fisher's fiducial distribution" [Neyman (1941)]. Its development has proceeded from Fisher (1930) though recent contributions, just to name a few, of Efron (1993, 1998), Fraser (1991, 1996), Lehmann (1993), Schweder and Hjort (2002) and others. There is renewed interest in CDs [Schweder and Hjort (2002)], partly because "statisticians will be asked to solve bigger and more complicated problems" [Efron (1998)] and the development of CDs might hold a key to "our profession's 250-year search for a dependable objective Bayes theory" [Efron (1998) and Schweder and Hjort (2002)].

This paper is mainly focused on some new developments on the "combination" aspect of CDs, where two natural approaches of combining CD information from independent studies are considered. The first approach is from the $p$-value combination scheme which dates back to Fisher (1932); see also Littell and Folks (1973) and Marden (1991), among many others. The second approach is analogous to multiplying likelihood functions in parametric inference. The two approaches are compared in the case of combining asymptotic normality based CDs. We require the resulting function of combined CDs to be a CD (or an asymptotic CD) so that it can be used later on to make inferences, store information or combine information in a sequential way.

For this purpose, we adopt a formal definition of CD developed by and presented in Schweder and Hjort (2002), and extend it to obtain a formal definition of asymptotic confidence distributions (aCDs). Suppose $X_1, X_2, \ldots, X_n$ are $n$ independent random draws from a population $F$ and $\mathcal{X}$ is the sample space corresponding to the data set $\mathbf{X}_n = (X_1, X_2, \ldots, X_n)^T$. Let $\theta$ be a parameter of interest associated with $F$ ($F$ may contain other nuisance parameters), and let $\Theta$ be the parameter space.

DEFINITION 1.1.   A function $H_n(\cdot) = H_n(\mathbf{X}_n, \cdot)$ on $\mathcal{X} \times \Theta \to [0, 1]$ is called a confidence distribution (CD) for a parameter $\theta$ if (i) for each given $\mathbf{X}_n \in \mathcal{X}$, $H_n(\cdot)$ is a continuous cumulative distribution function; (ii) at the true parameter value $\theta = \theta_0$, $H_n(\theta_0) = H_n(\mathbf{X}_n, \theta_0)$, as a function of the sample $\mathbf{X}_n$, has the uniform distribution $U(0, 1)$.

The function $H_n(\cdot)$ is called an asymptotic confidence distribution (aCD) if requirement (ii) above is replaced by (ii)′: at $\theta = \theta_0$, $H_n(\theta_0) \xrightarrow{W} U(0, 1)$ as $n \to +\infty$, and the continuity requirement on $H_n(\cdot)$ is dropped.

We call, when it exists, $h_n(\theta) = H_n'(\theta)$ a CD density, also known as a confidence density in the literature. This CD definition is the same as in Schweder and Hjort (2002), except that we suppress possible nuisance parameter(s) for notational simplicity. Our version, which was developed independently of Schweder and Hjort (2002), was motivated by the observation (1.1) below. For every $\alpha$ in $(0, 1)$, let $(-\infty, \xi_n(\alpha)]$ be a $100\alpha\%$ lower-side confidence interval, where $\xi_n(\alpha) = \xi_n(\mathbf{X}_n, \alpha)$ is continuous and increasing in $\alpha$ for each sample $\mathbf{X}_n$. Then $H_n(\cdot) = \xi_n^{-1}(\cdot)$ is a CD in the usual Fisherian sense. In this case,

$$\{\mathbf{X}_n : H_n(\theta) \leq \alpha\} = \{\mathbf{X}_n : \theta \leq \xi_n(\alpha)\}$$

(1.1)

$$\text{for any } \alpha \text{ in } (0,1) \text{ and } \theta \text{ in } \Theta \subseteq \mathbb{R}.$$

Thus, at $\theta = \theta_0$, $\Pr\{H_n(\theta_0) \leq \alpha\} = \alpha$ and $H_n(\theta_0)$ is $U(0,1)$ distributed. Definition 1.1 is very convenient for the purpose of verifying if a particular function is a CD or an aCD.

The notion of a CD (or aCD) is attractive for the purpose of combining information. The main reasons are that there is a wealth of information on $\theta$ inside a CD, the concept of CD (and particularly aCD) is quite broad, and the CDs are relatively easy to construct and interpret. Section 2 provides a brief review of materials related to the CDs along these views. See Schweder and Hjort (2002) for an expanded discussion of the concept of CDs and the information contained in CDs.

The main developments are in Sections 3 and 4. We provide in Section 3 a general recipe by adopting a general $p$-value combination scheme. Section 3.1 derives an optimal method for combining CDs associated with the same parameter, where the optimality is in terms of the Bahadur slope. The optimal scheme is notably different from that for combining $p$-values. Section 3.2 proposes adaptive combination methods, in the setting where the parameter values in some of the prior studies are not necessarily the same as the parameter value in the current study. The properties of adaptive consistency and adaptive efficiency are discussed. Analogous to combining likelihood functions in likelihood inference, we study in Section 4 a combination approach of multiplying CD densities. There we also provide a comparison of the two different CD-combining approaches in the case of normal type aCDs. Section 5 illustrates the methodology through three examples, each of which has individual significance. The proofs are in the Appendix.

**2. Examples and inferential information contained in a CD.** The notion of CDs and aCDs covers a broad range of examples, from regular parametric cases to $p$-value functions, normalized likelihood functions, bootstrap distributions and Bayesian posteriors, among others.

EXAMPLE 2.1. *Normal mean and variance.* Suppose $X_1, X_2, \ldots, X_n$ is a sample from $N(\mu, \sigma^2)$, with both $\mu$ and $\sigma^2$ unknown. A CD for $\mu$ is $H_n(y) = F_{t_{n-1}}(\frac{y - \overline{X}}{s_n/\sqrt{n}})$, where $\overline{X}$ and $s_n^2$ are, respectively, the sample mean and variance, and $F_{t_{n-1}}(\cdot)$ is the cumulative distribution function of the Student $t_{n-1}$-distribution. A CD for $\sigma^2$ is $H_n(y) = 1 - F_{\chi_{n-1}^2}(\frac{(n-1)s_n^2}{y})$ for $y \geq 0$, where $F_{\chi^2_{n-1}}(\cdot)$ is the cumulative distribution function of the $\chi^2_{n-1}$-distribution.

EXAMPLE 2.2. *p-value function.* For any given $\tilde{\theta}$, let $p_n(\tilde{\theta}) = p_n(\mathbf{X}_n, \tilde{\theta})$ be a *p*-value for a one-sided test $K_0: \theta \leq \tilde{\theta}$ versus $K_1: \theta > \tilde{\theta}$. Assume that the *p*-value is available for all $\tilde{\theta}$. The function $p_n(\cdot)$ is called a *p-value function* [Fraser (1991)]. Typically, at the true value $\theta = \theta_0$, $p_n(\theta_0)$ as a function of $\mathbf{X}_n$ is exactly (or asymptotically) $U(0,1)$-distributed. Also, $H_n(\cdot) = p_n(\cdot)$ for every fixed sample is almost always a cumulative distribution function. Thus, usually $p_n(\cdot)$ satisfies the requirements for a CD (or aCD).

EXAMPLE 2.3. *Likelihood functions.* There is a connection between the concepts of aCD and various types of likelihood functions, including likelihood functions in single parameter families, profile likelihood functions, Efron's implied likelihood function and Schweder and Hjort's reduced likelihood function, and so on. In fact, one can easily conclude from Theorems 1 and 2 of Efron (1993) that in an exponential family, both the profile likelihood and the implied likelihood [Efron (1993)] are aCD densities after a normalization. Singh, Xie and Strawderman (2001) provided a formal proof, with some specific conditions, which shows that $e^{\ell_n^*(\theta)}$ is proportional to an aCD density for the parameter $\theta$, where $\ell_n^*(\theta) = \ell_n(\theta) - \ell_n(\hat{\theta})$, $\ell_n(\theta)$ is the log-profile likelihood function, and $\hat{\theta}$ is the maximum likelihood estimator of $\theta$. Schweder and Hjort (2002) proposed the reduced likelihood function, which itself is proportional to a CD density for a specially transformed parameter. Also see Welch and Peers (1963) and Fisher (1973) for earlier accounts of likelihood function based CDs in single parameter families.

EXAMPLE 2.4. *Bootstrap distribution.* Let $\hat{\theta}$ be a consistent estimator of $\theta$. In the basic bootstrap methodology the distribution of $\hat{\theta} - \theta$ is estimated by the bootstrap distribution of $\hat{\theta}_B - \hat{\theta}$, where $\hat{\theta}_B$ is the estimator $\hat{\theta}$ computed on a bootstrap sample. An aCD for $\theta$ is $H_n(y) = P_B(\hat{\theta}_B \geq 2\hat{\theta} - y) = 1 - P_B(\hat{\theta}_B - \hat{\theta} \leq \hat{\theta} - y)$, where $P_B(\cdot)$ is the probability measure induced by bootstrapping. As $n \to \infty$, the limiting distribution of normalized $\hat{\theta}$ is often symmetric. In this case, due to the symmetry, the raw bootstrap distribution $H_n(y) = P_B(\hat{\theta}_B \leq y)$ is also an aCD for $\theta$.

Other examples include a second-order accurate CD of the population mean based on Hall's [Hall (1992)] second-order accurate transformed *t*-statistic, an aCD of the correlation coefficient based on Fisher's *z*-score function, among many others. See Schweder and Hjort (2002) for more examples and extended discussion.

A CD contains a wealth of information, somewhat comparable to, but different than, a Bayesian posterior distribution. A CD (or aCD) derived from a likelihood function can also be interpreted as an objective Bayesian posterior. We give a brief summary below of information in a CD related to some basic elements of inference. The reader can find more details in Singh, Xie and Strawderman ([2001](#)). This information is also scattered around in earlier publications, for example, in Fisher ([1973](#)), Fraser ([1991](#), [1996](#)) and Schweder and Hjort ([2002](#)), among others.

- *Confidence interval*. From the definition, it is evident that the intervals $(-\infty, H_n^{-1}(1-\alpha)]$, $[H_n^{-1}(\alpha), +\infty)$ and $(H_n^{-1}(\alpha/2), H_n^{-1}(1-\alpha/2))$ provide $100(1-\alpha)\%$-level confidence intervals of different kinds for $\theta$, for any $\alpha \in (0,1)$. The same is true for an aCD, where the confidence level is achieved in limit.
- *Point estimation*. Natural choices of point estimators of the parameter $\theta$, given $H_n(\theta)$, include the median $M_n = H_n^{-1}(1/2)$, the mean $\bar{\theta}_n = \int_{-\infty}^{+\infty} t \, dH_n(t)$ and the maximum point of the CD density $\hat{\theta}_n = \arg\max_\theta h_n(\theta)$, $h_n(\theta) = H_n'(\theta)$. Under some modest conditions one can prove that these point estimators are consistent plus more.
- *Hypothesis testing*. From a CD, one can obtain $p$-values for various hypothesis testing problems. Fraser ([1991](#)) developed some results on such a topic through $p$-value functions. The natural line of thinking is to measure the support that $H_n(\cdot)$ lends to a null hypothesis $K_0 : \theta \in C$. We perceive two types of support: 1. *Strong-support* $p_s(C) = \int_C dH_n(\theta)$. 2. *Weak-support* $p_w(C) = \sup_{\theta \in C} 2\min(H_n(\theta), 1 - H_n(\theta))$. If $K_0$ is of the type $(-\infty, \theta_0]$ or $[\theta_0, \infty)$ or a union of finitely many intervals, the strong-support $p_s(C)$ leads to the classical $p$-values. If $K_0$ is a singleton, that is, $K_0$ is $\theta = \theta_0$, then the weak-support $p_w(C)$ leads to the classical $p$-values.

## 3. Combination of CDs through a monotonic function.

In this section we consider a basic methodology for combining CDs which essentially originates from combining $p$-values. However, there are some new twists, modifications and extensions. Here one assumes that some past studies (with reasonably sensible results) on the current parameter of interest exist. The CDs to be combined may be based on different models. A nice feature of this combination method is that, after combination, the resulting function is always an exact CD if the input CDs from the individual studies are exact. Also, it does not require any information regarding how the input CDs were obtained. Section 3.1 considers the perfect situation when the common parameter had the same value in all previous studies on which the CDs are based. Section 3.2 presents an adaptive combination approach which works asymptotically, even when there exist some "wrong CDs" (CDs with underlying true parameter values different from $\theta_0$). For clarity, the presentation in this section is restricted to CDs only. The entire development holds for aCDs with little or no modification.

3.1. *CD combination and Bahadur efficiency.* Let $H_1(y), \ldots, H_L(y)$ be $L$ independent CDs, with the same true parameter value $\theta_0$ (sample sizes are suppressed in the CD notation in the rest of this paper). Suppose $g_c(u_1, \ldots, u_L)$ is any continuous function from $[0,1]^L$ to $\mathbb{R}$ that is monotonic in each coordinate. A general way of combining, depending on $g_c(u_1, \ldots, u_L)$, can be described as follows: Define $H_c(u_1, \ldots, u_L) = G_c(g_c(u_1, \ldots, u_L))$, where $G_c(\cdot)$ is the continuous cumulative distribution function of $g_c(U_1, \ldots, U_L)$, and $U_1, \ldots, U_L$ are independent $U(0,1)$ distributed random variables. Denote

$$(3.1) \qquad H_c(y) = H_c(H_1(y), \ldots, H_L(y)).$$

It is easy to verify that $H_c(y)$ is a CD function for the parameter $\theta$. We call $H_c(y)$ a combined CD. If the objective is only to get a combined aCD, one may also allow the above $g_c$ function to involve sample estimates.

Let $F_0(\cdot)$ be any continuous cumulative distribution function and $F_0^{-1}(\cdot)$ be its inverse function. A convenient special case of the function $g_c$ is

$$(3.2) \qquad g_c(u_1, u_2, \ldots, u_L) = F_0^{-1}(u_1) + F_0^{-1}(u_2) + \cdots + F_0^{-1}(u_L).$$

In this case, $G_c(\cdot) = F_0 * \cdots * F_0(\cdot)$, where $*$ stands for convolution. Just like the $p$-value combination approach, this general CD combination recipe is simple and easy to implement. Some examples of $F_0$ are:

- $F_0(t) = \Phi(t)$ is the cumulative distribution function of the standard normal. In this case

$$H_{NM}(y) = \Phi\left(\frac{1}{\sqrt{L}}[\Phi^{-1}(H_1(y)) + \Phi^{-1}(H_2(y)) + \cdots + \Phi^{-1}(H_L(y))]\right).$$

- $F_0(t) = 1 - e^{-t}$, for $t \geq 0$, is the cumulative distribution function of the standard exponential distribution (with mean 1). Or, $F_0(t) = e^t$, for $t \leq 0$, which is the cumulative distribution function of the mirror image of the standard exponential distribution. In these cases the combined CDs are, respectively,

$$H_{E1}(y) = P\left(\chi^2_{2L} \leq -2 \sum_{i=1}^{L} \log(1 - H_i(y))\right)$$

and

$$H_{E2}(y) = P\left(\chi^2_{2L} \geq -2 \sum_{i=1}^{L} \log H_i(y)\right),$$

where $\chi^2_{2L}$ is a $\chi^2$-distributed random variable with $2L$ degrees of freedom. The recipe for $H_{E2}(y)$ corresponds to Fisher's recipe of combining $p$-values [Fisher (1932)].

- $F_0(t) = \frac{1}{2}e^t \mathbb{1}_{(t \leq 0)} + (1 - \frac{1}{2}e^{-t})\mathbb{1}_{(t \geq 0)}$, denoted as $DE(t)$ from now on, is the cumulative distribution function of the standard double exponential distribution. Here $\mathbb{1}_{(\cdot)}$ is the indicator function. In this case the combined CD is

$$H_{DE}(y) = DE_L(DE^{-1}(H_1(y)) + \cdots + DE^{-1}(H_L(y))),$$

where $DE_L(t) = DE * \cdots * DE(t)$ is the convolution of $L$ copies of $DE(t)$.

Lemma 3.1 next gives an iterative formula to compute $DE_L(t)$. One critical fact of this lemma is that the exponential parts of the tails of $DE_L(t)$ are the same as those of $DE(t)$. The proof of Lemma 3.1 is in the Appendix.

LEMMA 3.1. *For $t > 0$ we have*

$$1 - DE_L(t) = DE_L(-t) = \frac{1}{2}V_L(t)e^{-t},$$

*where $V_L(t)$ is a polynomial of degree $L - 1$. This sequence of polynomials satisfies the following recursive relation: for $k = 2, 3, \ldots, L$,*

$$2V_k(t) = V_{k-1}(t) + \int_0^t [V_{k-1}(s) - V'_{k-1}(s)]\,ds$$

$$+ \int_0^\infty [V_{k-1}(s) + V_{k-1}(t+s) - V'_{k-1}(s)]e^{-2s}\,ds.$$

*In particular, $V_1(t) = 1$ and $V_2(t) = 1 + t/2$.*

Littell and Folks ([1973](#)) established an optimality property, in terms of Bahadur slope, within the class of combined $p$-values based on monotonic combining functions. Along the same line, we establish below an optimality result for the combination of CDs.

Following Littell and Folks ([1973](#)), we define the concept of Bahadur slope for a CD:

DEFINITION 3.1. Let $n$ be the sample size corresponding to a CD function $H(\cdot)$. We call a nonnegative function $S(t) = S(t; \theta_0)$ the Bahadur slope for the CD function $H(\cdot)$ if for any $\varepsilon > 0$, $S(-\varepsilon) = -\lim_{n \to +\infty} \frac{1}{n} \log H(\theta_0 - \varepsilon)$ and $S(\varepsilon) = -\lim_{n \to +\infty} \frac{1}{n} \log\{1 - H(\theta_0 + \varepsilon)\}$ almost surely.

The Bahadur slope gives the rate, in exponential scale, at which $H(\theta_0 - \varepsilon)$ and $1 - H(\theta_0 + \varepsilon)$ go to zero. The larger the slope, the faster its tails decay to zero. In this sense, a CD with a larger Bahadur slope is asymptotically more efficient as a "distribution-estimator" for $\theta_0$.

Suppose $n_1, n_2, \ldots, n_L$, the sample sizes behind the CDs $H_1(y), H_2(y), \ldots, H_L(y)$, go to infinity at the same rate. For notational simplicity, replace $n_1$ by $n$ and write $n_j = \{\lambda_j + o(1)\}n$ for $j = 1, 2, 3, \ldots, L$; we always have $\lambda_1 = 1$. Let $S_j(t)$ be the Bahadur slope for $H_j(y)$, $j = 1, 2, \ldots, L$, and $S_c(t)$ be the Bahadur slope for their combined CD, say $H_c(y)$. The next theorem provides an upper bound for $S_c(t)$ [i.e., the fastest possible decay rate of $H_c(y)$ tails] and indicates when it is achieved. Its proof can be found in the Appendix.

THEOREM 3.2.   *Under $\theta = \theta_0$, for any $\varepsilon > 0$, as $n \to +\infty$,*

$$-\liminf \frac{1}{n} \log H_c(\theta_0 - \varepsilon) \le \sum_{j=1}^{L} \lambda_j S_j(-\varepsilon)$$

*and*

$$-\liminf \frac{1}{n} \log(1 - H_c(\theta_0 + \varepsilon)) \le \sum_{j=1}^{L} \lambda_j S_j(\varepsilon)$$

*almost surely. If the slope function $S_c(\cdot)$ exists,*

$$S_c(-\varepsilon) \le \sum_{j=1}^{L} \lambda_j S_j(-\varepsilon) \Big/ \sum_{j=1}^{L} \lambda_j \quad \text{and} \quad S_c(\varepsilon) \le \sum_{j=1}^{L} \lambda_j S_j(\varepsilon) \Big/ \sum_{j=1}^{L} \lambda_j.$$

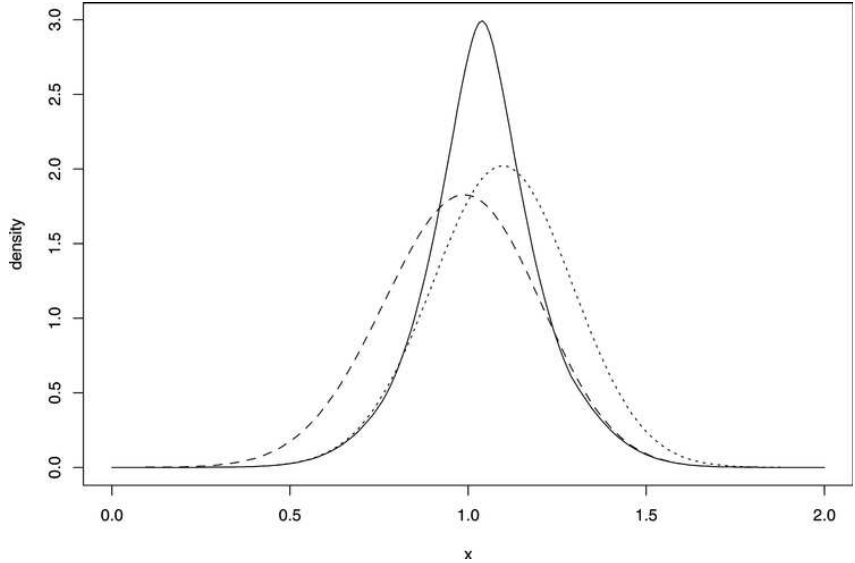

FIG. 1.   *A typical figure of density plots of $H_{DE}(\cdot)$, $H_{E1}(\cdot)$ and $H_{E2}(\cdot)$, when we combine independent CDs for the common mean parameter $\mu$ of $N(\mu, 1.0)$ and $N(\mu, 1.5^2)$; the true $\mu = 1$ and the sample sizes are $n_1 = 30$ and $n_2 = 40$. The solid, dotted and dashed curves are the density curves of $H_{DE}(\cdot)$, $H_{E1}(\cdot)$ and $H_{E2}(\cdot)$, respectively.*

*Furthermore, almost surely, for any $\varepsilon > 0$,*

$$S_{E1}(\varepsilon) = \sum_{j=1}^{L} \lambda_j S_j(\varepsilon) \Big/ \sum_{j=1}^{L} \lambda_j \quad and \quad S_{E2}(-\varepsilon) = \sum_{j=1}^{L} \lambda_j S_j(-\varepsilon) \Big/ \sum_{j=1}^{L} \lambda_j;$$

$$S_{DE}(-\varepsilon) = \sum_{j=1}^{L} \lambda_j S_j(-\varepsilon) \Big/ \sum_{j=1}^{L} \lambda_j \quad and \quad S_{DE}(\varepsilon) = \sum_{j=1}^{L} \lambda_j S_j(\varepsilon) \Big/ \sum_{j=1}^{L} \lambda_j.$$

*Here, $S_{E1}(t)$, $S_{E2}(t)$ and $S_{DE}(t)$ are the Bahadur slope functions of combined CDs $H_{E1}(x)$, $H_{E2}(x)$ and $H_{DE}(x)$, respectively.*

This theorem states that the $DE(t)$ based combining approach is, in fact, optimal in terms of achieving the largest possible value of Bahadur slopes on both sides. The two combined CDs $H_{E1}(y)$ and $H_{E2}(y)$ can achieve the largest possible slope value only in one of the two regions, $\theta > \theta_0$ or $\theta < \theta_0$. This phenomenon is illustrated in Figure 1.

Note that, in the *p*-value combination case, Littell and Folks ([1973](#)) established that Fisher's way of combination is optimal in terms of having the largest possible Bahadur slope. To our knowledge, no one has considered the $DE(t)$ based combination rule when combining *p*-values. There is a notable difference between combining *p*-values and CDs. While for CDs one cares about the decay rates of both tails, separately, a typical *p*-value concept either involves only one tail of the distribution of a test statistic or lumps its two tails together. The $DE(t)$ based combination rule is quite natural when combining CDs, but not when combining *p*-values.

### 3.2. *Adaptive combination.*

The development in Section 3.1 is under the assumption that all the CDs $H_1(y), \ldots, H_L(y)$ are for the same parameter $\theta$, with identical true value $\theta = \theta_0$. There may be doubts about the validity of this assumption. For instance, let $H_1(y)$ be a CD for $\theta$ with true value $\theta_0$, based on a current study of sample size $n_1$. Let $H_2(y), \ldots, H_L(y)$ be available CDs on $\theta$ based on previous (independent) studies involving sample sizes $\{n_2, \ldots, n_L\}$. One could be less than certain that all earlier values of $\theta$ were indeed equal to the current value, $\theta = \theta_0$. It will be problematic if one combines all the available CDs when some of the studies had the underlying value $\theta \neq \theta_0$. Indeed, the resulting function of combination will not even be an aCD (under the true value $\theta = \theta_0$). This can be demonstrated by a simple example of combining two CDs: $H_1(y) = \Phi(\frac{y - \hat{\theta}_1}{\sigma/\sqrt{n}})$ and $H_2(y) = \Phi(\frac{y - \hat{\theta}_2}{\sigma/\sqrt{n}})$, where $\sigma$ is known, $\sqrt{n}(\hat{\theta}_1 - \theta_0) \sim N(0, \sigma^2)$, and $\sqrt{n}(\hat{\theta}_2 - \tilde{\theta}_0) \sim N(0, \sigma^2)$, for some $\tilde{\theta}_0 \neq \theta_0$. The combined outcome by the normal-based approach, $H_{NM}(y) = \Phi(\frac{2y - \hat{\theta}_1 - \hat{\theta}_2}{\sqrt{2}\sigma/\sqrt{n}})$, is not uniformly distributed, even in limit ($n \to \infty$) when $y = \theta_0$. In this section we propose adaptive combination approaches to remedy this problem and go on to establish related optimality properties under a large sample setting.

Let $H_1(\cdot)$ be a CD with true parameter value $\theta = \theta_0$, based on a current study. The $L-1$ CDs from the previous studies on $\theta$ are separated into two sets:

$$\mathcal{H}_0 = \{H_j : H_j \text{ has the underlying value of } \theta = \theta_0, j = 2, \ldots, L\},$$

$$\mathcal{H}_1 = \{H_j : H_j \text{ has the underlying value of } \theta \neq \theta_0, j = 2, \ldots, L\}.$$

The set $\mathcal{H}_0$ contains the "right" CDs and $\mathcal{H}_1$ contains the "wrong" CDs. We assume, however, the information about $\mathcal{H}_0$ and $\mathcal{H}_1$ is unavailable to us.

The development of a general adaptive combination recipe starts with an extension of the general combination method of Section 3.1, which includes a set of weights $\boldsymbol{\omega} = (\omega_1, \ldots, \omega_L)$, $\omega_1 \equiv 1$. Our intention is to select a set of adaptive weights that can filter out the "wrong" CDs (in $\mathcal{H}_1$) and keep the "right" CDs (in $\mathcal{H}_0$) asymptotically.

Although it could be much more general, for simplicity, we let

$$(3.3) \qquad g_{c,\omega}(u_1, \ldots, u_L) = \sum_{j=1}^{L} \omega_j F_0^{-1}(u_j),$$

where $F_0(\cdot)$ is a continuous cumulative distribution function with a bounded density. Let $G_{c,\omega}(t) = F_0(t) * F_0(\frac{t}{\omega_2}) * \cdots * F_0(\frac{t}{\omega_L})$ and $H_{c,\omega}(u_1, \ldots, u_L) = G_{c,\omega}(g_{c,\omega}(u_1, \ldots, u_L))$. We define

$$H_{c,\omega}(y) = H_{c,\omega}(H_1(y), H_2(y), \ldots, H_L(y)).$$

Define the weight vector $\boldsymbol{\omega}_0$ as $(1, \omega_2^{(0)}, \ldots, \omega_L^{(0)})$, where $\omega_j^{(0)} = 1$ for $H_j \in \mathcal{H}_0$, and $\omega_j^{(0)} = 0$ for $H_j \in \mathcal{H}_1$. The combined CD function $H_c^{(0)}(y) = H_{c,\omega_0}(y)$ is our target which combines $H_1$ with the CDs in $\mathcal{H}_0$. Of course, we lack the knowledge of $\mathcal{H}_0$ and $\mathcal{H}_1$, so $\boldsymbol{\omega}_0$ is unknown. Thus, we need to determine the adaptive weights, denoted by $\boldsymbol{\omega}^* = (1, \omega_2^*, \ldots, \omega_L^*)$, converging to $\boldsymbol{\omega}_0$, from the available information in the CDs. Let $H_c^*(y) = H_{c,\omega^*}(y)$. One would hope that $H_c^*(y)$ is at least an aCD.

DEFINITION 3.2.   A combination method is adaptively consistent if $H_c^*(y)$ is an aCD for $\theta = \theta_0$.

Suppose $n_1, n_2, \ldots, n_L$ go to infinity at the same rate. Again, we let $n = n_1$ and write $n_i = \{\lambda_i + o(1)\}n$, for $i = 1, 2, \ldots, L$; $\lambda_1 = 1$. We define below adaptive slope efficiency.

DEFINITION 3.3. A combination method is adaptively slope efficient if for any $\varepsilon > 0$,

$$-\lim_{n \to +\infty} \frac{1}{n} \log H_c^*(-\varepsilon) = \sum_{j \,:\, H_j \in \{H_1\} \cup \mathcal{H}_0} \lambda_j S_j(-\varepsilon),$$

$$-\lim_{n \to +\infty} \frac{1}{n} \log\{1 - H_c^*(\varepsilon)\} = \sum_{j \,:\, H_j \in \{H_1\} \cup \mathcal{H}_0} \lambda_j S_j(\varepsilon).$$

Here $S_j(t)$ is the Bahadur slope of $H_j(x)$, all assumed to exist.

Let $I_i$ be a confidence interval derived from $H_i(\cdot)$, for $i = 1, 2, \ldots, L$. Suppose, as $n \to \infty$, the lengths of these $L$ confidence intervals all go to zero almost surely. Then, for all large $n$, it is expected that $I_1 \cap I_j = \varnothing$ for $j$ such that $H_j \in \mathcal{H}_1$. This suggests that in order to get rid of the CDs in $\mathcal{H}_1$ when $n$ is large, we should take $\omega_j^* = \mathbb{1}_{(I_1 \cap I_j \neq \varnothing)}$ for $j = 2, \ldots, L$, where $\mathbb{1}_{(\cdot)}$ is the indicator function. With this choice of data-dependent weights, we have the following theorem. The theorem can be easily proved using the Borel–Cantelli lemma and its proof is omitted.

THEOREM 3.3. *Assume that the intervals $I_j$ lie within an $\varepsilon$-neighborhood of the corresponding true value of $\theta$, for all large $n$ almost surely, and for any fixed $\varepsilon > 0$. In addition, assume that, for each $j$ such that $H_j \in \mathcal{H}_0$, we have*

(3.4) $$\sum_{n=1}^{+\infty} P(I_1 \cap I_j = \varnothing) < +\infty.$$

*Then if $\omega_j^* = I_{(I_1 \cap I_j \neq \varnothing)}$ for $j = 1, 2, \ldots, L$, we have $\sup_y |H_c^*(y) - H_c^{(0)}(y)| = 0$, for all large $n$ almost surely.*

Note that $H_c^{(0)}(y)$ is the "target" combined CD. From Theorem 3.3 we immediately have the following corollary.

COROLLARY 3.4. *Under the assumptions of Theorem 3.3, the adaptive combination recipe described in this section, with the adapting weights $\omega_j^* = \mathbb{1}_{(I_1 \cap I_j \neq \varnothing)}$ for $j = 2, \ldots, L$, is adaptively consistent. Furthermore, if $F_0(t) = DE(t)$ in (3.3), this combination method is also adaptively slope efficient.*

REMARK 3.1. A simple example is to take $I_i = (H_i^{-1}(\alpha_n/2), H_i^{-1}(1 - \alpha_n/2))$, $i = 1, 2, \ldots, L$. It follows that $P(I_1 \cap I_j = \varnothing) \leq 1 - (1 - \alpha_n)^2 \leq 2\alpha_n$ for each $j$ such that $H_j \in \mathcal{H}_0$. Thus, $\sum_{n=1}^{+\infty} \alpha_n < \infty$ is a sufficient condition for (3.4). However, this bound is typically very conservative. To see this, consider the basic example of $z$-based CD for unknown normal mean with known variance $\sigma^2$. Let $H_i(y) = \Phi(\sqrt{n}(y - \overline{X}_i)/\sigma)$, where $\overline{X}_i$ is the sample mean in the $i$th study, and $z_{\alpha_n/2}$ is the normal critical value of level $\alpha_n/2$. We have $P(I_1 \cap I_j = \varnothing) = 2(1 - \Phi(\sqrt{2}z_{\alpha_n/2}))$, which could be a lot smaller than $2\alpha_n$. Considering this issue, we recommend a somewhat substantial value for $\alpha_n$ in applications.

A feature that may be regarded as undesirable with the above adaptive method is the fact that it assigns weights either 0 or 1 to the CDs. We propose below the use of kernel function based weights, which take values between 0 and 1. Under some regularity conditions, we will show that the weighted adaptive combination method with the kernel based weights is adaptively consistent and "locally efficient" (Theorem 3.5).

Let $K(t)$ be a symmetric kernel function, $\int K(t)\,dt = 1$, $\int tK(t)\,dt = 0$ and $\int t^2 K(t)\,dt = 1$. In the present context we also require that the tails of the kernel function tend to zero at an exponential rate. Some examples are the normal kernel $K(t) = \phi(t)$, the triangle kernel and the rectangular kernel function, among others.

In order to use the kernel function, some measure of "distance" between $H_1(y)$ and $H_j(y)$, $j = 2, \ldots, L$, is needed. For illustrative purposes, we use $\hat{\theta}_1 - \hat{\theta}_j$, where $\hat{\theta}_i$, $i = 1, \ldots, L$, are point estimators obtained from $H_i(y)$, respectively. We assume $\hat{\theta}_i$, for $i = 1, 2, \ldots, L$, converge in probability to their respective underlying values of $\theta$, say $\theta_{0,i}$, at the same polynomial rate. For $i = 1$ or $i$ such that $H_i \in \mathcal{H}_0$, $\theta_{0,i} \equiv \theta_0$. Let $b_n \to 0$ be such that $|\hat{\theta}_i - \theta_{0,i}| = o_p(b_n)$. We define the kernel function based weights as

$$(3.5) \qquad \omega_j^* = K\left(\frac{\hat{\theta}_1 - \hat{\theta}_j}{b_n}\right) \Big/ K(0) \qquad \text{for } j = 1, 2, \ldots, L.$$

Among many other possibilities, one set of convenient choices is $\hat{\theta}_i = H_i^{-1}(\frac{1}{2})$ and $b_n = \sqrt{R_n}$, where $R_n = H_1^{-1}(\frac{3}{4}) - H_1^{-1}(\frac{1}{4})$ is the interquartile range of $H_1(y)$.

Under the above setting, we have the following theorem; its proof is in the Appendix.

THEOREM 3.5. Let $\delta_n > 0$ be a sequence such that $H_i(\theta_0 \pm \delta_n)$ are bounded away from 0 and 1, in probability, for $i = 1$ and $i$ with $H_i \in \mathcal{H}_0$. Suppose $F_0$ in (3.3) is selected such that $\min\{F_0(t), 1 - F_0(t)\}$ tends to zero, exponentially fast as $|t| \to \infty$. Then, with $\omega_j^*$ as defined in (3.5), one has

$$\sup_{x \in [\theta_0 - \delta_n, \theta_0 + \delta_n]} |H_c^*(y) - H_c^{(0)}(y)| \to 0 \qquad \text{in probability, as } n \to +\infty.$$

Theorem 3.5 suggests that in a local neighborhood of $\theta_0$, $H_c^*(y)$ and $H_c^{(0)}(y)$ are close for large $n$. Recall that $H_c^{(0)}(y)$ is the target that combines $H_1$ with the CDs in $\mathcal{H}_0$. The following conclusion is immediate from Theorem 3.5.

COROLLARY 3.6. *Under the setting of Theorem 3.5, with $\omega_j^*$ as in* (3.5), *the adaptive combination method described in this section is adaptively consistent.*

The result in Theorem 3.5 is a local result depending on $\delta_n$, which is typically $O(n^{-1/2})$. For a set of general kernel weights of the form (3.5), we cannot get an anticipated adaptive slope efficiency result for the adaptive DE combination method. But, for the rectangular kernel, this optimality result does hold, since in this case the weight $\omega_j^*$ becomes either 1 or 0 for all large $n$, almost surely. The proof of the following corollary is similar to that of Theorem 3.4 and is omitted.

COROLLARY 3.7. *Under the setting of Theorem 3.5, with $\omega_j^*$ as in* (3.5) *and $K(t) = \{1/(2\sqrt{3})\}\mathbb{1}_{(|t|<\sqrt{3})}$, the adaptive combination method described in this section with $F_0(t) = DE(t)$ is adaptively slope efficient if $\sum_{n=1}^{\infty} P(|\hat{\theta}_j - \theta_{0,j}| > \frac{\sqrt{3}}{2}b_n) < \infty$ for $j = 1, 2, \ldots, L$.*

## 4. Combination of CDs through multiplying CD densities.

Normalized likelihood functions (as a function of the parameter) are an important source of obtaining CD or aCD densities. In fact, it was Fisher who prescribed the use of normalized likelihood functions for obtaining his fiducial distributions; see, for example, Fisher (1973). Multiplying likelihood functions from independent sources constitutes a standard method for combining parametric information. Naturally, this suggests multiplying CD densities and normalizing to possibly derive combined CDs as follows:

$$(4.1) \qquad H_P(\theta) = \int_{(-\infty,\theta)\cap\Theta} h^*(y)\,dy \Big/ \int_{(-\infty,\infty)\cap\Theta} h^*(y)\,dy,$$

where $h^*(y) = \prod_{i=1}^{L} h_i(y)$ and $h_i(y)$ are CD densities from $L$ independent studies. Schweder and Hjort (2002) suggested multiplying their reduced likelihood functions for combined estimation, which is closely related to the approach of (4.1). However, they did not require normalization and, strictly speaking, the reduced likelihood function, in general, is not a CD density for $\theta$ (it is only proportional to a CD density for a specially transformed parameter).

Unfortunately, the combined function $H_P(\cdot)$ may not necessarily be a CD or even an aCD function in general. But we do have some quite general affirmative results. We first present here a basic result pertaining to $H_P(\cdot)$. Let $T_1, T_2, \ldots, T_L$ be a set of statistics from $L$ independent samples. Suppose $\widetilde{H}_i(\cdot)$, $i = 1, 2, \ldots, L$, are the cumulative distribution functions of $T_i - \theta$ with density functions $\tilde{h}_i(\cdot)$ which are entirely free of parameters. Thus, one has $L$ CDs of $\theta$, given by $H_i(\theta) = 1 - \widetilde{H}_i(T_i - \theta)$ with corresponding CD densities $h_i(\theta) = \tilde{h}_i(T_i - \theta)$.

THEOREM 4.1.   *In the above setting, $H_P(\theta)$ is an exact CD of $\theta$.*

An elementary proof of this theorem is given in the Appendix. This theorem can also be proved using the general theory relating best equivariant procedures and Bayes procedures relative to right invariant Haar measure [see, e.g., Berger (1985), Chapter 6 for the Bayes invariance theory]. Using this Bayes-equivariance connection, or directly, one can also obtain an exact CD for the scale parameter $\theta$, but it requires one to replace $h^*(y)$ in (4.1) with $h^{**}(y)/y$, where $h^{**}(y) = \prod_{i=1}^{L}\{yh_i(y)\}$.

The original method of (4.1) does not yield an exact CD for a scale parameter. Let us consider a simple example.

EXAMPLE 4.1.   Consider the $U[0, \theta]$ distribution with unknown $\theta$. Let $H_i(\theta) = 1 - (\frac{Y_i}{\theta})^{n_i}$ over $\theta \geq Y_i$, $i = 1, 2$, be the input CDs, where $Y_1$ and $Y_2$ are maxima of two independent samples of sizes $n_1$ and $n_2$. The multiplication method (4.1) yields $H_P(\theta) = (\frac{Y}{\theta})^{n_1+n_2+1}$, over $\theta \geq Y = \max(Y_1, Y_2)$. This $H_P(\theta)$ is not an exact CD, though it is an aCD.

The setting for Theorem 4.1 is limited. But it allows an asymptotic extension that covers a wide range of problems, including those involving the normal and "heavy tailed" asymptotics, as well as other nonstandard asymptotics such as that in Example 4.1.

Let $\widetilde{H}_{i,a}$ be an asymptotic (weak limit) cumulative distribution function of $\xi_i = n_i^{\alpha}\frac{T_i - \theta}{V_i}$, where $(T_i, V_i)$ are statistics based on independent samples of sizes $n_i$, $i = 1, \ldots, L$. Denote $\tilde{h}_{i,a}(\cdot) = \widetilde{H}'_{i,a}(\cdot)$. One has aCD densities given by

$$(4.2) \qquad h_{i,a}(\theta) = \frac{n_i^{\alpha}}{V_i}\tilde{h}_{i,a}\left(n_i^{\alpha}\frac{T_i - \theta}{V_i}\right).$$

Let $\xi_i$ have uniformly bounded exact densities $\tilde{h}_{i,e}(\cdot)$, for $i = 1, \ldots, L$, and define $h_{i,e}(\cdot)$ as in (4.2) with $\tilde{h}_{i,a}$ replaced by $\tilde{h}_{i,e}(\cdot)$. Assume the regularity conditions: (a) $\tilde{h}_{i,e}(\cdot) \to \tilde{h}_{i,a}(\cdot)$ uniformly on compact sets. (b) $\tilde{h}_{i,e}(\cdot)$ are uniformly integrable. (c) $V_i \to \tau_i^2$, a positive quantity, in probability, for $i = 1, \ldots, L$. Define $H_P(\cdot)$ by (4.1) where $h_i(\cdot)$ is either $h_{i,a}(\cdot)$ or $h_{i,e}(\cdot)$.

THEOREM 4.2. *In the above setting, $H_P(\cdot)$ is an aCD.*

The proof is based on standard asymptotics using Theorem 4.1 on combination of $H_{i,a}$'s and is not presented here. We would like to remark that, due to the special form of the normal density, in the usual case of normal asymptotics ($\alpha = \frac{1}{2}$, $\widetilde{H}_{i,a} = \Phi$), the combined function $H_P(\cdot)$, with $h_i(\cdot) = h_{i,a}(\cdot)$ in (4.1), is an aCD without requiring the regularity conditions (a) and (b).

For the purpose of comparing the two different combination approaches given by (3.1) and (4.1), we now specialize to asymptotic normality based aCDs where both methods can apply. Let $\xi_i = \sqrt{n_i} \frac{T_i - \theta}{V_i}$. The normality based aCD is $H_{i,a}(\theta) = 1 - \widetilde{H}_{i,a}(\xi_i) = 1 - \Phi(\xi_i)$ with aCD density $h_{i,a}(\theta) = \tilde{h}_{i,a}(\xi_i) = \frac{\sqrt{n_i}}{V_i} \phi(\xi_i)$. Consider the combined function $H_P(\cdot)$ with input aCD densities $h_{i,a}(\cdot)$ or $h_{i,e}(\cdot)$. It is straightforward to verify that $H_P(\cdot)$ in this special case is the same as (or asymptotically equivalent to)

$$H_{AN}(\theta) = 1 - \Phi\left(\left[\sum_{i=1}^{L} \frac{n_i}{V_i}\right]^{1/2} (\hat{\theta}_c - \theta)\right),$$

where $\hat{\theta}_c = \left(\sum_{i=1}^{L} \frac{n_i}{V_i} T_i\right) / \sum_{i=1}^{L} \left(\frac{n_i}{V_i}\right)$ is the asymptotically optimal linear combination of $T_i$, $i = 1, 2, \ldots, L$. In light of this remark, it is evident that the large sample comparison presented below between $H_{AN}$ and $H_{DE}$ also holds between $H_P$ and $H_{DE}$. Note that $H_{AN}$ is, in fact, a member of the rich class of combining methods introduced in Section 3.1, where we pick $g_c(u_1, \ldots, u_L) = \sum_{i=1}^{L} [(\frac{n_i}{V_i})^{1/2} \Phi^{-1}(u_i)]$ in (3.1).

The concept of the Bahadur slope, which is at the heart of Section 3, is still well defined for aCDs and $H_{DE}(\cdot)$ still has the slope optimality. However, the concept of slope loses its appeal on aCDs since one can alter the slope of an aCD by tampering with its tails, while keeping it an aCD. Nevertheless, if the input CDs are the normal based aCDs mentioned above, it is straightforward to show that $H_{DE}$ and $H_{AN}$ have the same slope (achieving the upper bound). This result is noteworthy in the special case when $T_i$'s are means of normal samples and $V_i = \sigma_i^2$ are the known variances, where $H_{AN}$ is an exact CD. In this case, $H_{AN}$ is derived from a UMVUE estimator.

Next we address the following question: How do $H_{DE}$ and $H_{AN}$ compare in terms of the lengths of their derived equal-tail, two-sided asymptotic confidence intervals? Let us define $\ell_{DE}(\alpha) = H_{DE}^{-1}(1-\alpha) - H_{DE}^{-1}(\alpha)$ and $\ell_{AN}(\alpha) = H_{AN}^{-1}(1-\alpha) - H_{AN}^{-1}(\alpha)$. Also, we assume, as in Section 3, that $n_1 = n$ and $n_j/n$ are bounded below and above, for $j = 2, \ldots, L$. Let $\lim^*$ be the limit as $(n, \alpha) \to (\infty, 0)$. A key aspect of $\lim^*$ is the fact that it allows $\alpha$ to converge to 0, at an arbitrary rate; of course, slow rates are better from a practical viewpoint. The proof is given in the Appendix.

THEOREM 4.3. $\lim^*[\ell_{DE}(\alpha)/\ell_{AN}(\alpha)] = 1$ *in probability.*

Hence, for combining large sample normality based aCDs, the three combining methods $H_{DE}$, $H_P$ and $H_{AN}$ are equivalent in the above sense. When the input aCDs are derived from profile likelihood functions, the $H_P$-method amounts to multiplying profile likelihood functions, in which case (in view of the standard likelihood inference) $H_P$ may have a minor edge over $H_{DE}$ when a finer comparison is employed. On the other hand, $H_{DE}$ has its global appeal, especially when nothing is known about how the input CDs were derived. It always returns an exact CD when the input CDs are exact. Also, $H_{DE}$ preserves the second-order accuracy when the input CDs are second-order accurate (a somewhat nontrivial result, not presented here). Aspects of second-order asymptotics on $H_P$ are not known to us, while $H_{AN}$ ignores second-order corrections.

The adaptive combining in Section 3.2 carries over to $H_{AN}$, since $H_{AN}(\cdot)$ is a member of the rich class of combining methods introduced there. Also, one can turn $H_P$ into an adaptively combined CD by replacing $h^*(y)$ in (4.1) with $h_\omega^*(y)$, where $h_\omega^*(y) = \prod_{i=1}^{L} h_i^{\omega_i}(y)$ or $\prod_{i=1}^{L} h_i(\omega_i y)$. The adaptive weights $\omega_i$ are chosen such that $\omega_i \to 1$ for the "right" CDs (in $\mathcal{H}_0$) and $\omega_i \to 0$ for the "wrong" CDs (in $\mathcal{H}_1$). Some results along the line of Section 3.2 can be derived.

We close this section with the emerging recommendation that while normal type aCDs can be combined by any of the methods $H_{DE}$, $H_P$ or $H_{AN}$, exact CDs and higher-order accurate CDs should generally be combined by the $DE$ method.

## 5. Examples.

5.1. *The common mean problem.* The so-called common mean problem of making inference on the common mean, say $\mu$, of two or more normal populations of possibly different variances, also known as the Behrens–Fisher problem, has attracted a lot of attention in the literature. In the large sample setting, it is well known that the Graybill–Deal estimator, $\hat{\mu}_{GD} = \{(n_2/s_2^2)\overline{X}_1 + (n_1/s_1^2)\overline{X}_2\}/\{(n_1/s_1^2) + (n_2/s_2^2)\}$, is asymptotically efficient. In the small sample setting, there is still research going on attempting to find efficient exact confidence intervals for $\mu$. In particular, Jordan and Krishnamoorthy [[1996](), through combining statistics] and Yu, Sun and Sinha [[1999](), through combining two-sided $p$-values] proposed efficient exact confidence intervals for the common mean $\mu$; however, there is a small but nonzero chance that these intervals do not exist.

Let us consider the CD based method, first under large sample settings. In this case, we start with normal based aCDs $H_{1a}(y) = \Phi(\frac{y - \overline{X}_1}{s_1/\sqrt{n_1}})$ and $H_{2a}(y) = \Phi(\frac{y - \overline{X}_2}{s_2/\sqrt{n_2}})$. Following Section 4.2, we would like to prescribe the combined CDs $H_{AN}(y)$ [or $H_P(y)$, which is the same]. It is interesting to note that this combined CD is the same as the CD directly derived from the Graybill–Deal estimator. Thus, the confidence intervals derived from $H_{AN}(\theta)$ are asymptotically shortest.

If one wants to obtain exact confidence intervals for $\mu$, one can turn to the recipe prescribed in Section 3.1. Clearly, exact CDs for $\mu$ based on two independent normal samples are $H_1(y) = F_{t_{n_1-1}}(\frac{y-\overline{X}_1}{s_1/\sqrt{n_1}})$ and $H_2(y) = F_{t_{n_2-1}}(\frac{y-\overline{X}_2}{s_2/\sqrt{n_2}})$, respectively; see Example 2.1. By Theorem 3.2, the DE based approach will be Bahadur optimal among all exact CD based approaches of Section 3. The resulting exact confidence interval for $\mu$, with coverage $1-\alpha$, is $(q_{\alpha/2}, q_{1-\alpha/2})$, where $q_s$ is the $s$-quantile of the CD function $H_{DE}(y)$. This exact confidence interval for $\mu$ always exists at every level $\alpha$.

We carried out a simulation study of 1000 replications to examine the coverage of the CD based approaches, under three sets of sample sizes $(n_1, n_2) = (3, 4), (30, 40)$ or $(100, 140)$ and two sets of (true) variances $(\sigma_1^2, \sigma_2^2) = (1, 1.5^2)$ or $(1, 3.5^2)$. The coverage of constructed 95% confidence intervals is right on target around 95% in the six cases for the $H_{DE}$ based exact method. However, the Graybill–Deal (i.e., aCD $H_{AN}$ or $H_P$) based method leads to serious under-coverage (84.8% and 85.9%) in the two cases with small sample sizes $(n_1, n_2) = (3, 4)$, and notable under-coverage (93.3% and 93.6%) in the two cases with moderate sample sizes $(n_1, n_2) = (30, 40)$. So, in small sample cases, the exact CD based approach is substantially better, in terms of coverage.

Theorem 4.3 suggests that, under a large sample setting, the DE based approach and the Graybill–Deal estimator (equivalently, $H_{AN}$ or $H_P$ based approach will have similar lengths for confidence intervals with high asymptotic coverage. We carried out a simulation study to compare the lengths in the two cases with large sample size $(n_1, n_2) = (100, 140)$, at confidence level 95%. We found that the lengths corresponding to the $H_{DE}$ based method, on average, are slightly higher than those corresponding to the Graybill–Deal estimator, but they are not too far apart. The average ratio of the lengths, in the 1000 simulations, is 1.034 for $(\sigma_1^2, \sigma_2^2) = (1, 1.5^2)$ and 1.081 for $(\sigma_1^2, \sigma_2^2) = (1, 3.5^2)$. Similar ratios were also obtained for the 90% and 99% confidence intervals under the same setting. The simulation results seem to endorse our recommendation at the end of Section 4.

5.2. *Adaptive combination of odds ratios in ulcer data.* Efron [[1993](#), Section 5] studied an example of combining independent studies. The example concerns a randomized clinical trial studying a new surgical treatment for stomach ulcers [Kernohan, Anderson, McKelvey and Kennedy ([1984](#))] in which there are 9 successes and 12 failures among 21 patients in treatment groups, and 7 successes and 17 failures among 24 patients in the control group. The parameter of interest is the log odds ratio of the treatment.

Based on the data, the estimated log odds ratio is $\hat{\theta}_1 = \log(\frac{9}{12}/\frac{7}{17}) = 0.600$, with estimated standard error $\hat{\sigma}_1 = (\frac{1}{9} + \frac{1}{12} + \frac{1}{7} + \frac{1}{17})^{1/2} = 0.629$. In addition to Kernohan's trial, there were 40 other randomized trials of the same treatment between 1980 and 1989 [see Table 1 of Efron (1996) for the complete data]. The question of interest is how to combine the information in these 40 studies with that in Kernohan's trial. Efron (1993) employed an empirical Bayes approach, where he used a Bayes rule to combine the implied likelihood function of Kernohan's trial $L_x^*(\theta) \approx \frac{1}{\hat{\sigma}_1}\phi(\frac{\theta-\hat{\theta}_1}{\hat{\sigma}_1})$ with a prior distribution $\pi_e(\theta) \propto \sum_{j=2}^{41} \frac{1}{\hat{\sigma}_i}\phi(\frac{\theta-\hat{\theta}_j}{\hat{\sigma}_j})$. Here $\phi(t)$ is the density function of the standard normal distribution, and $\hat{\theta}_j$ and $\hat{\sigma}_j$, $j = 2, \ldots, 41$, are the estimators of the log odds ratios and standard errors in the 40 other clinical trials. To obtain meaningful estimates of $\hat{\theta}_j$ and $\hat{\sigma}_j$ in the analysis, nine entries of zero were changed to 0.5; see Efron (1993).

We re-study this example, utilizing the purely frequentist CD combination approach. Under the standard assumption that the data in each of these 41 independent clinical trials are from a four-category multinomial distribution, it is easy to verify that $H_j(y) = \Phi(\frac{y-\hat{\theta}_j}{\hat{\sigma}_j})$, $j = 1, 2, \ldots, 41$, are a set of first-order normal aCDs of the 41 clinical trials. We use the combined aCD $H_{AN}$ (i.e., taking $g_c(u_1, \ldots, u_L) = \sum_{i=1}^{41}[\frac{1}{\hat{\sigma}_i}\Phi^{-1}(u_i)]$ in (3.1)), both with and without adaptive weights, to summarize the combined information. Although there is no way to theoretically compare our approach with Efron's empirical Bayes approach, we will discuss the similarities and differences of the final outcomes from these two alternative approaches.

First, let us temporarily assume that the underlying values of $\theta$ in these 41 clinical trials are all the same. So, each trial receives the same weight in combination. In this case, the combined aCD is $H_{AN}^S(\theta) = \Phi(\{\sum_{i=1}^{41} \frac{1}{\hat{\sigma}_i}\frac{\theta-\hat{\theta}_i}{\hat{\sigma}_i}\}/\{\sum_{i=1}^{41}\frac{1}{\hat{\sigma}_i^2}\}^{1/2}) = \Phi(7.965(\theta + 0.8876))$. The density curve of $H_{AN}^S(\theta)$ is plotted in Figure 2(a), along with the posterior density curve (dashed line) obtained from Efron's empirical Bayes approach. For easy illustration, we also include (in each plot) two dotted curves that correspond to the aCD density of Kernohan's trial $h_1(\theta) = H_1'(\theta)$ and the average aCD densities of the previous 40 trials $g_a(\theta) = \frac{1}{41}\sum_{i=2}^{41}\frac{1}{\hat{\sigma}_i}\phi(\frac{\theta-\hat{\theta}_i}{\hat{\sigma}_i})$; note that $h_1(\theta) \approx L_x^*(\theta)$, Efron's (1993) implied likelihood $L_x^*(\theta)$, and $g_a(\theta) \propto \pi_e(\theta)$, the empirical prior used in Efron (1993). It is clear in Figure 2(a) that the aCD curve of $H_{AN}^S(\theta)$ is too far to the left, indicating a lot of weight has been given to the 40 other trials. We believe that the assumption of the same underlying values of $\theta$ in all of these 41 clinical trials is too strong; see also Efron (1996).

A more reasonable assumption is that some of the 40 other trials may not have the same underlying true $\theta$ as in Kernohan's trial. It is sensible to use the adaptive combination methods proposed in Section 3.2, which downweight or exclude the trials with the underlying parameter value away from that of Kernohan's trial. Three sets of adaptive weights are considered: one set of normal kernel weights $\omega_j^N \propto \phi(\frac{M_1 - M_j}{\sqrt{R}})$, and two sets of 0 or 1 adaptive weights $\omega_j^I = \mathbb{1}_{(I_1 \cap I_j \neq \varnothing)}$ with $I_i = (H_i^{-1}(\alpha_n/2), H_i^{-1}(1 - \alpha_n/2))$. Here $M_i = H_i(\frac{1}{2})$ is the median of the aCD $H_i$ of the $i$th trial, and following Remark 3.1, we take $\alpha_n = 0.25$ and 0.30, respectively, in the two sets of $\omega_j^I$'s. The three corresponding combined CDs are, respectively, $H_{AN}^N(\theta) = \Phi(5.7788(\theta - 0.1029))$, $H_{AN}^{I_1}(\theta) = \Phi(5.4007(\theta - 0.1199))$ and $H_{AN}^{I_2}(\theta) = \Phi(5.3051(\theta - 0.1698))$. Their density curves are plotted in Figure 2(b)–(d). We also tried triangle and rectangular kernel weights, but the results were very similar and are not presented here. A noteworthy feature of the combined CD density curves is the following: when all the weights are equal, the combined curve puts little mass to the right of 0, while all the rest put substantial mass to the right of 0.

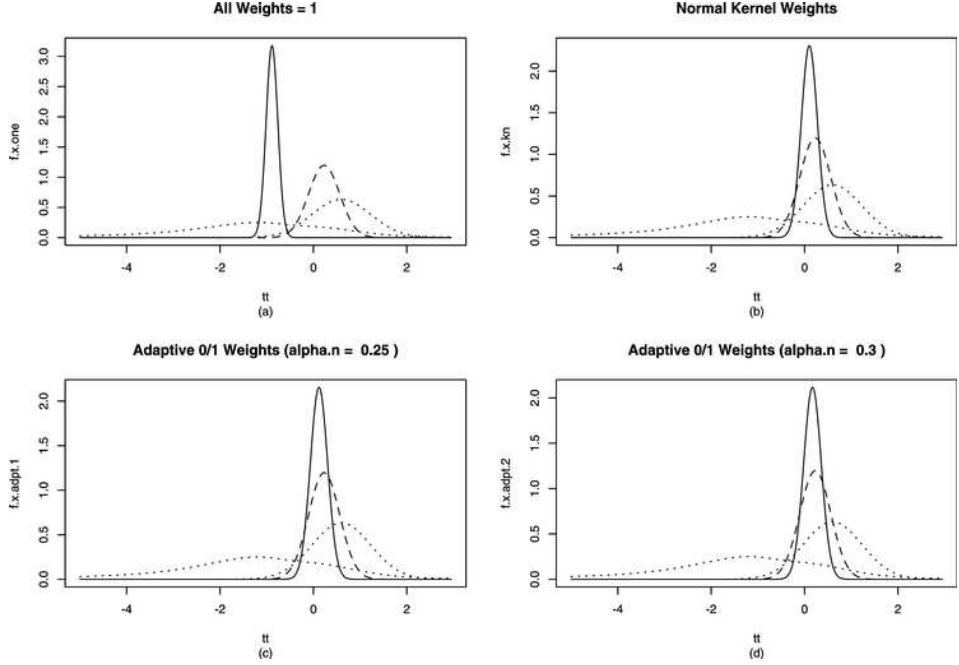

Fig. 2. *The solid curves in* (a)–(d) *are combined CD density curves, combined with* (a) *equal weight to all 41 trials,* (b) *normal kernel weights,* (c) *0–1 weights with* $\alpha_n = 0.25$, (d) *0–1 weights with* $\alpha_n = 0.30$. *The dashed curves are the posterior density function* (*approximated*) *from Figure 4 of Efron* (1993). *The two dotted curves* (*with peaks from right to left*) *correspond to the aCD density of Kernohan's trial* $h_1(\theta)$ [*i.e., Efron's* (1993) *implied likelihood* $L_x^*(\theta)$] *and the average aCD densities of the 40 other trials* [ *proportional to the empirical prior* $\pi_e(\theta)$ *used in Efron* (1993)].

Comparing the three adaptively combined CDs with the posterior distribution obtained by Efron [[1993](#), Figure 4] on the same data set [Figure 2(b)–(d)], we find that they all have very similar features. Their density curves all peak within a small range and between the peaks of $h_1(\theta)$ and $g_a(\theta)$, actually much closer to $h_1(\theta)$, reflecting intuitive intent of such combinations. But there is also a quite notable difference. The spans of all the combined CD densities are smaller than that of the posterior density function. Note that in Efron's empirical Bayes approach, all 40 other trials have equal contributions (as a part of the prior) to the final posterior through Bayes formula. In the adaptive approach, the trials closer to Kernohan's trial have more contribution (i.e., higher weights) than those trials farther away. It seems that much more information from the 40 other clinical trials, especially those with $H_j(\theta)$ closer to $H_1(\theta)$, has been drawn in the adaptive CD combination method.

### 5.3. Computationally intense methodology on a large data set.

One can utilize CDs to find a way to apply statistical methodology involving heavy computations on a large data set. Here, we illustrate the "split and combine" approach. We divide the data into smaller data sets; after analyzing each sub-data set separately, we can piece together useful information through combining CDs. For a computationally intense methodology, such a method can result in tremendous saving. Suppose the number of steps involved in a statistical methodology is $cn^{1+a}$, $n$ being the size of the data set, $a > 0$. Suppose the data set is divided into $k$ pieces, each of size $\frac{n}{k}$. The number of steps involved in carrying out the method on each subset is $c(\frac{n}{k})^{1+a}$. Thus, the total number of steps is $ck(\frac{n}{k})^{1+a} = \frac{cn^{1+a}}{k^a}$. If the effort involved in combining CDs is ignored, there is a saving by a factor of $k^a$. We think that the information loss due to this approach will be minimal. One question is how to divide the data. Simple approaches include dividing the data based on their indices (time or natural order index), random sampling or some other natural groupings.

For the purpose of demonstration, let us consider a $U$-statistic based robust multivariate scale proposed by Oja ([1983](#)). Let $\{X_1, \ldots, X_n\}$ be a two-dimensional data set. Oja's robust multivariate scale is $\mathcal{S}_n = \text{median}\{\text{areas of all } \binom{n}{3}$ triangles formed by 3 data points$\}$. For any given three data points $(x_1, y_1)$, $(x_2, y_2)$ and $(x_3, y_3)$, the area of their triangle is given by $\frac{1}{2}|\det(\mathbf{t}_1 \mathbf{t}_2 \mathbf{t}_3)|$, where $\mathbf{t}_l = (1, x_l, y_l)'$, for $l = 1, 2, 3$. To make inference on this scale parameter, it is natural to use the bootstrap. But obtaining the bootstrap density of $\mathcal{S}_n$ is a formidable task when the sample size $n$ is large. For example, even with $n = 48$ the computation of $\mathcal{S}_n$ involves evaluating the area of $48 \times 47 \times 46/6 = 17{,}296$ triangles. With 5000 (say) repeated bootstrapping, the total number of triangle areas needed to be evaluated is 86.48 million. If one adopts the "split and combine" approach discussed above, say, randomly splitting the data set of size 48 into two data sets of size 24, the number of triangle areas needed to be evaluated is $2 \times 5000 \times (24 \times 23 \times 22/6) = 20.24$ million. This is less than $1/4$ of the 86.48 million. If we randomly split the data set of size 48 into three data sets of size 16 each, the number of triangle areas needed to be evaluated is $3 \times 5000 \times (16 \times 15 \times 14/6) = 8.4$ million, less than $1/10$ of the 86.48 million. Since bootstrap density functions are aCDs, the bootstrap density functions obtained from each sub-dataset can be combined together, using the techniques of combining CDs. The combined CD can be used to make inferences on the robust multivariate scale $\mathcal{S}_n$.

Figure 3 plots bootstrap density functions of the robust multivariate scale $\mathcal{S}_n$ based on a simulated two-dimensional data set of size 48. The data set was generated with the $s$th observation being $(z_s^{[1]} + z_s^{[2]}, z_s^{[1]} - z_s^{[2]})$, where $z_s^{[1]}$ and $z_s^{[2]}$, $s = 1, \ldots, 48$, are simulated from the Cauchy distributions with parameters center $= 0$ and scale $= 1$ and center $= 1$ and scale $= 1.3$, respectively. The solid curve in Figure 3(a) is the bootstrap density function of $\mathcal{S}_n$ based on 5000 bootstrapped samples. It took 67720.75 seconds to generate the density curve on our Ultra 2 Sparc Station using Splus. The dotted and broken curves are the bootstrap density functions of $\mathcal{S}_n$ using the "split and combine" method, where we split the sample randomly into two and three pieces, respectively. It took 12,647.73 and 5561.19 seconds to generate these two density plots, including the combining part. Hence, the "split and combine" method used less than 1/4 and 1/10 of the time, respectively! From Figure 3, it is apparent that all three curves are quite alike. They all seem to capture essential features of the bootstrap distribution of the robust multivariate scale $\mathcal{S}_n$.

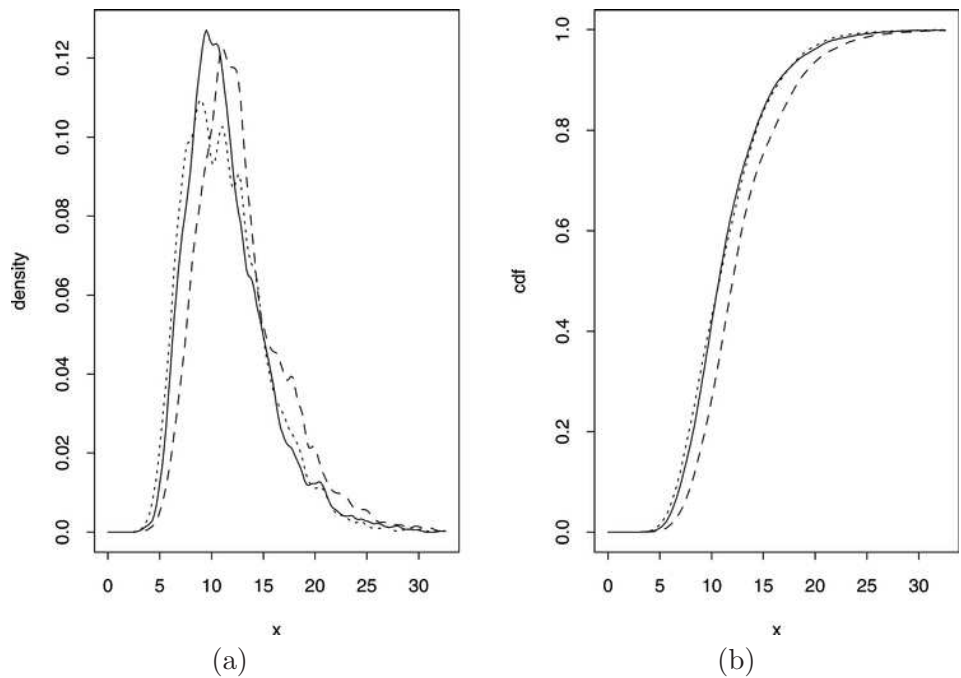

(a)  (b)

FIG. 3. *Figure* (a) *is for density curves and* (b) *is for cumulative distribution curves. The solid, dotted and dashed-line curves correspond to methods described in the main context with no split, split into two pieces and split into three pieces, respectively.*

## APPENDIX

### A.1. Proofs in Section 3.

PROOF OF LEMMA 3.1. Let $T$ and $W$ be independent random variables, such that $T$ has the standard double exponential distribution while $W$ satisfies, for $t > 0$, $P(W \leq -t) = P(W > t) = \frac{1}{2} V_k(t)e^{-t}$, where $V_k(\cdot)$ is a polynomial of finite degree. We write, for $t > 0$, $P(T + W > t) = \frac{1}{4}P(T + W > t | T > 0, W > 0) + \frac{1}{4}P(T + W > t | T > 0, W < 0) + \frac{1}{4}P(T + W > t | T < 0, W > 0) = I + II + III$ (say). Now,

$$I = \frac{1}{4}P(W > t | T > 0, W > 0) + \frac{1}{4}P(T + W > t, W \leq t | T > 0, W > 0)$$

$$= \frac{1}{4}V_k(t)e^{-t} + \frac{1}{4}\int_0^t e^{-(t-s)}[V_k(s) - V_k'(s)]e^{-s}\,ds$$

$$= \frac{1}{4}e^{-t}\left[V_k(t) + \int_0^t [V_k(s) - V_k'(s)]\,ds\right],$$

$$II = \frac{1}{4}P(T > t - W | T > 0, W < 0)$$

$$= \frac{1}{4}P(T > t + W | T > 0, W > 0)$$

$$= \frac{1}{4}\int_0^\infty e^{-(t+s)}[V_k(s) - V_k'(s)]e^{-s}\,ds$$

$$= \frac{1}{4}e^{-t}\left[\int_0^\infty \{V_k(s) - V_k'(s)\}e^{-2s}\,ds\right].$$

Similarly, $III = \frac{1}{4}e^{-t}[\int_0^\infty V_k(s+t)e^{-2s}\,ds]$. The proof is concluded by letting $W$ have the distribution of $\sum_1^k Y_i$, for $k = 1, 2, \ldots$, successively. $\quad\square$

Before we prove Theorem 3.2, we first prove a lemma. This lemma borrows an idea from Littell and Folks (1973).

LEMMA A.1. *Let* $H_c(u_1, \ldots, u_L)$ *be a function from* $[0, 1]^L$ *to* $[0, 1]$, *monotonically nondecreasing in each coordinate. Also, suppose* $H_c(U_1, \ldots, U_L)$ *has the* $U(0, 1)$ *distribution when* $U_1, \ldots, U_L$ *are independent* $U(0, 1)$ *random* *variables. Then for any* $u_1, \ldots, u_L$ *in* $[0, 1]$, $H_c(u_1, \ldots, u_L) \geq \prod_{\ell=1}^{L} u_\ell$ *and* $1 - H_c(u_1, \ldots, u_L) \geq \prod_{\ell=1}^{L}(1 - u_l)$.

PROOF. In view of the monotonicity, it follows that $\{U_1 \leq u_1, U_2 \leq u_2, \ldots, U_L \leq u_L\}$ implies $\{H_c(U_1, \ldots, U_L) \leq H_c(u_1, \ldots, u_L)\}$. The first claim follows if we take $U_1, U_2, \ldots, U_L$ as independent $U(0, 1)$ random variables. The second claim can be proved similarly via the facts that $\{U_1 \geq u_1, \ldots, U_L \geq u_L\}$ implies $\{1 - H_c(U_1, \ldots, U_L) \leq 1 - H_c(u_1, \ldots, u_L)\}$ and that $1 - H_c(U_1, \ldots, U_L)$ follows $U(0, 1)$ distribution when $U_1, \ldots, U_L$ are independent $U(0, 1)$ random variables. $\square$

PROOF OF THEOREM 3.2. In Lemma A.1 take $u_1 = H_1(x), \ldots, u_L = H_L(x)$. The first result follows immediately from Lemma A.1. Note $m/n_1 \to \sum \lambda_j$, where $m = n_1 + n_2 + \cdots + n_L$. The second result follows from the first result.

The next two equalities related to $H_{E1}$ and $H_{E2}$ can be obtained from direct calculations, appealing to the fact that the upper tail area of the $\chi^2_{2L}$-distribution satisfies $\lim_{y \to +\infty} \frac{1}{y} \log P(\chi^2_{2L} > y) = -\frac{1}{2}$, where $\chi^2_{2L}$ is a $\chi^2_{2L}$-distributed random variable. Note, by Lemma 3.1, $\lim_{y \to +\infty} \frac{1}{y} \log DE_L(-y) = \lim_{y \to +\infty} \frac{1}{y} \log(1 - DE_L(y)) = -1$. Using this, it is seen that the last two claims also hold. $\square$

The proof of Theorem 3.5 critically depends on the following lemma.

LEMMA A.2. *Under the condition of Theorem 3.5, with* $\omega_n^*$ *as in* (3.5), *one has*

$$\sup_t |G_{c,\omega^*}(t) - G_{c,\omega^{(0)}}(t)| \to 0 \qquad \text{in probability,}$$

*where* $\boldsymbol{\omega}^{(0)} = (1, \omega_2^{(0)}, \ldots, \omega_L^{(0)})$, $\omega_j^{(0)} = 1$ *if* $H_j \in \mathcal{H}_0$ *and* $\omega_j^{(0)} = 0$ *if* $H_j \in \mathcal{H}_1$.

PROOF. Let $\mathbf{Y} = \{Y_1, \ldots, Y_L\}$ be i.i.d. r.v.'s having the distribution $F_0$, independent of the original data. Clearly $\omega_i^* \to \omega_i^{(0)}$ in probability, for $i = 2, \ldots, L$. Note that when $|\omega_i^* - \omega_i^{(0)}| < \delta$, for a $\delta > 0$ and all $i = 2, \ldots, L$, we have

$$P_{\mathbf{Y}}\left(Y_1 + \sum_{i=2}^{L} \omega_i^{(0)} Y_i + \delta \sum_{i=2}^{L} |Y_i| \leq t\right) \leq P_{\mathbf{Y}}\left(Y_1 + \sum_{i=2}^{L} \omega_i^* Y_i \leq t\right)$$

$$\leq P_{\mathbf{Y}}\left(Y_1 + \sum_{i=2}^{L} \omega_i^{(0)} Y_i - \delta \sum_{i=2}^{L} |Y_i| \leq t\right).$$

Using standard arguments, one deduces that, as $|\delta| \to 0$,

$$\sup_t \left| P_{\mathbf{Y}}\left(Y_1 + \sum_{i=2}^{L} \omega_i^{(0)} Y_i \leq t\right) - P_{\mathbf{Y}}\left(Y_1 + \sum_{i=2}^{L} \omega_i^{(0)} Y_i + \delta \sum_{i=2}^{L} |Y_i| \leq t\right) \right| \to 0.$$

The lemma follows from combining the above assertions. $\square$

PROOF OF THEOREM 3.5. For $H_i \in \mathcal{H}_0$, $\sup_{|x-\theta_0| \leq \delta_n} |\omega_i^* - \omega_i^{(0)}| \times F_0^{-1}(H_i(x)) \overset{p}{\to} 0$ using the condition on $\delta_n$. For $H_i \in \mathcal{H}_1$ and $|x - \theta_0| \leq \delta_n$, $\min\{H_i(x), 1 - H_i(x)\}$ tends to 0 at an exponential rate. Therefore, $F_0^{-1}(H_i(x)) = O(n)$, since the tails of $F_0$ decay at an exponential rate as well. Using the assumed condition on $b_n$ and the kernel function, we deduce that $\omega_i^* \to 0$ in probability, faster than any polynomial rate, for $i$ such that $H_i \in \mathcal{H}_1$. Thus, for such an $i$, $\sup_{|y-\theta_0| \leq \delta_n} |\omega_i^* F_0^{-1}(H_i(y))| \to 0$ in probability. The theorem now follows, utilizing Lemma A.2. $\square$

**A.2. Proofs in Section 4.**

PROOF OF THEOREM 4.1. For simplicity, we assume $\Theta = (-\infty, +\infty)$; other situations can be dealt with similarly. Note we only need to prove that $H_P(\theta_0)$ is $U(0, 1)$ distributed.

Define an $L-1$ random vector $\mathbf{Z} = (Z_1, Z_2, \ldots, Z_{L-1})^T$, where $Z_j = T_j - T_L$, for $j = 1, 2, \ldots, L-1$. So the joint density function of $\mathbf{Z}$ is $f_Z(\mathbf{z}) = f_Z(z_1, z_2, \ldots, z_{L-1}) = \int_{-\infty}^{+\infty} \prod_{j=1}^{L-1} \tilde{h}_j(z_j + u) \tilde{h}_L(u) \, du$ and the conditional density of $T_L - \theta_0$, given $\mathbf{Z}$, is $f_{T_L|Z}(t) = \prod_{j=1}^{L-1} \tilde{h}_j(Z_j + t) \tilde{h}_L(t) / f_Z(\mathbf{Z})$. Also, for each given $\mathbf{Z}$, we define a decreasing function $K_{\mathbf{Z}}(\gamma) = \int_{\gamma}^{+\infty} \prod_{j=1}^{L-1} \tilde{h}_j(Z_j + u) \tilde{h}_L(u) \, du$. It is clear that

$$H_P(\theta) = K_{\mathbf{Z}}(T_L - \theta) / f_Z(\mathbf{Z}).$$

So for any $s$, $0 < s < 1$, we have

$$\begin{aligned}
P\{H_P(\theta_0) \leq s\} &= P\{T_L - \theta_0 \geq K_{\mathbf{Z}}^{-1}(s f_Z(\mathbf{Z}))\} \\
&= E[P\{T_L - \theta_0 \geq K_{\mathbf{Z}}^{-1}(s f_Z(\mathbf{Z})) | \mathbf{Z}\}] \\
&= E\left[\int_{K_{\mathbf{Z}}^{-1}(s f_Z(\mathbf{Z}))}^{\infty} \frac{\prod_{j=1}^{L-1} \tilde{h}_j(Z_j + t) \tilde{h}_L(t)}{f_Z(\mathbf{Z})} \, dt\right] \\
&= E\left[\int_0^{s f_Z(Z)} \frac{1}{f_Z(\mathbf{Z})} \, du\right] = s,
\end{aligned}$$

where the fourth equality is due to a monotonic variable transformation in the integration: $u = K_Z(t)$. $\square$

PROOF OF THEOREM 4.3. Let $\theta_0$ be the common value of the parameter. To prove the claim, we show that for any $\varepsilon > 0$, $P(H_{DE}(\theta_0 + (1+\varepsilon)\ell_{AN}(\alpha)) > 1 - \alpha) \to 1$, $P(H_{DE}(\theta_0 + (1-\varepsilon)\ell_{AN}(\alpha)) < 1 - \alpha) \to 1$, and similar results on the lower side of $\theta_0$. We establish the first claim below; others can be proved similarly.

Let us note that, under $\lim^*$,

$$\sum_{i=1}^{L} DE^{-1}\left(\Phi\left(\frac{\sqrt{n_i}}{\tau_i}[\theta_0 + (1+\varepsilon)\ell_{AN}(\alpha) - T_i]\right)\right)$$

$$= \sum_{i=1}^{L} DE^{-1}\left(\Phi\left(\frac{\sqrt{n_i}}{\tau_i}[(1+\varepsilon)\ell_{AN}(\alpha)][1 + o_p(1)]\right)\right)$$

$$= \sum_{i=1}^{L} \frac{1}{2}\{\sqrt{n_i}[(1+\varepsilon)\ell_{AN}(\alpha)][1 + o_p(1)]/\tau_i\}^2 = [(1+\varepsilon)^2 z_\alpha^2/2][1 + o_p(1)].$$

Thus, by Lemma 3.1,

$$1 - H_{DE}(\theta_0 + (1+\varepsilon)\ell_{AN}(\alpha)) = o_p(\alpha^{1+\varepsilon}). \qquad \square$$

**Acknowledgments.** The comments of the referees and the editors of *The Annals of Statistics* greatly helped us improve the focus, presentation and some contents of the article.

K. SINGH
M. XIE
W. E. STRAWDERMAN
DEPARTMENT OF STATISTICS
RUTGERS—THE STATE UNIVERSITY
    OF NEW JERSEY
HILL CENTER, BUSCH CAMPUS
PISCATAWAY, NEW JERSEY 08854
USA
E-MAIL: kesar@stat.rutgers.edu
E-MAIL: mxie@stat.rutgers.edu
E-MAIL: straw@stat.rutgers.edu